\documentclass[showpacs,preprintnumbers,amsmath,amssymb,12pt]{article}
\usepackage{graphicx,amsfonts,amssymb,amsthm,amsmath,graphics,epsfig,pstricks,fancyhdr,fancybox}
\usepackage{dcolumn}
\usepackage{bm}

\newtheorem{thm}{Theorem}[section]

\newtheorem{lem}[thm]{Lemma}
\newtheorem{prop}[thm]{Proposition}

\newcommand{\refeq}[1]{~(\ref{#1})}
\newcommand{\myref}[1]{~\ref{#1}}
\newcommand{\mycite}[1]{~\cite{#1}}
\newcommand{\boxit}[1]{\vbox{\hrule\hbox{\vrule\kern5pt
     \vbox{\kern5pt#1\kern3pt}\kern5pt\vrule}\hrule}}

\newcommand{\sqr}[2]{{\vcenter{\vbox{\hrule height.#2pt
    \hbox{\vrule width.#2pt height#1pt \kern#1pt
    \vrule width.#2pt}\hrule height.#2pt}}}}

\newcommand{\att}[1]{\E\!\left[{#1}\right]}
\newcommand{\var}[1]{\V\left[{#1}\right]}
\newcommand{\pr}[1]{\P\left\{{#1}\right\}}

\newcommand{\indim}{\noindent{\bf Proof:}\hspace{0.2cm}}
\newcommand{\findim}{\hfill$\blacksquare$\vspace{0.3cm}\noindent}

\def\RE{{\mathbb R}}

\def\P{{\bm P}}
\def\Q{{\bf Q}}

\def\E{{\bm E}\,}
\def\V{{\bm V}}

\def\OM{\Omega}
\def\om{\omega}

\def\tok{\buildrel k \over \longrightarrow}

\def\va{\emph{rv}}

\def\FD{\emph{cdf}}

\def\mmax{k}
\def\Q0{\mathbb{Q}_0}

\begin{document}
\thispagestyle{empty}

\title{\LARGE \textbf{Taking rational numbers at random}}
\author{\textsc{Nicola Cufaro Petroni}\\
Department of Mathematics and \textsl{TIRES}, University of Bari\\
\textsl{INFN} Sezione di Bari\\
via E. Orabona 4, 70125 Bari, Italy\\
email: \textit{cufaro@ba.infn.it}}

\date{}

\maketitle

\begin{abstract}
\noindent We outline some simple prescriptions to define a
distribution on the set $\Q0$ of all the rational numbers in
$[0,1]$, and we then explore both a few properties of these
distributions, and the possibility of making these rational numbers
asymptotically equiprobable in a suitable sense. In particular it
will be shown that in the said limit -- albeit no uniform
distribution can be properly defined on $\Q0$ -- the probability
allotted to a single $q\in\Q0$ asymptotically vanishes, while that
of the subset of $\Q0$ falling in an interval $[a,b]$ goes to $b-a$.
We finally give some hints to completely sequencing without
repetitions the numbers in $\Q0$ as a prerequisite to the laying
down of more distributions on it
\end{abstract}

{
%\small\noindent MSC:

%\noindent PACS:

}

\section{Introduction}\label{intro}

What could the locution \emph{taking at random} possibly mean? In
its most general sense this would indicate that the drawing of an
outcome $\om$ out of a set $\OM$ is made according to any arbitrary
(but legitimate) probability measure assigned on the subsets of
$\OM$, and then that the usual precepts of the probability theory
are followed with the result that different probabilities are
normally allocated to distinct subsets of $\OM$. Traditionally
however the meaning of the said locution is more circumscribed and
stands rather for assuming that there is no reason to think that
there are preferred outcomes $\om\in\OM$, these being supposed
instead to be equally likely. This will be the meaning that we will
be interested in all along this paper, or -- whether this notion
will not be exactly applicable -- an asymptotic version of it in
some acceptable limiting sense

It is well known indeed that for the sets of real numbers our kind
of \emph{randomness} is enforced either by sheer
\emph{equiprobability} (on the finite sets), or by distribution
\emph{uniformity} (on the bounded, Lebesgue measurable, uncountable
sets). On the other hand infinite, countable sets and unbounded,
uncountable sets are both excluded from these egalitarian
probability attributions because their elements can be made neither
equiprobable (with a non vanishing probability), nor uniformly
distributed (with a non vanishing probability density). In these
occasions it is advisable instead to start with some proper (neither
equiprobable, nor uniform) probability distribution, and then to
inquire if and how this can be made ever closer -- in a suitable,
approximate sense -- either to an equiprobable or to an uniform one:
we will then respectively speak of \emph{asymptotic equiprobability}
and \emph{asymptotic uniformity}

The focus of our present inquiry, as will be elucidated in the
Section\myref{ratio}, are the rational numbers that -- even in a
bounded interval -- constitute an infinite, countable set, with a
few relevant, additional peculiarities due to their being also
everywhere dense among the real numbers. In the
Section\myref{ratiodistr} we will then supply a procedure to
attribute non vanishing probabilities to every rational number
$q=\,\!^n/_m$ in the interval $[0,1]$. The Section\myref{equinum} is
instead devoted to the aftermaths of supposing conditionally
equiprobable numerators $n$, and then the Section\myref{equi} will
show under which hypotheses our distributions can give rise to an
asymptotic equiprobability of the rationals in $[0,1]$ such that --
without pretending to have a uniform distribution on $\Q0$ -- the
probability allotted to a single $q\in\Q0$ vanishes in the limit,
while that of the subset of $\Q0$ falling in an interval $[a,b]$
goes to $b-a$. Several examples of denominator distributions are
elaborated in the Section\myref{denom} giving rise to a few closed
formulas, and finally in the Section\myref{concl} some concluding
remarks are added with a glimpse on the open problem of sequencing
all the rational numbers in $\Q0$

\section{Probability on rational numbers}\label{ratio}

Rational numbers are famously countable, and hence they can be put
in a sequence. Since however they are a dense subset of the real
numbers, every rational number is a cluster point, and hence no
sequence encompassing all of them can ever converge, not to say be
monotone. In any case their countability certainly allows the
allotment of discrete distributions with non vanishing probabilities
for every rational number: since they are infinite, however, they
can never be exactly \emph{equiprobable}. We will outline in the
forthcoming sections a simple procedure to give distributions on the
rationals in $[0,1]$, a set that we will shortly denote as
$\Q0=\mathbb{Q}\cap[0,1]$, and we will investigate if and how they
can be considered \emph{asymptotically equiprobable}. We will
refrain instead for the time being from extending our considerations
to the whole of $\mathbb{Q}$ only because in our opinion this -- at
the present stage of the inquiry -- would not add particular
insights to our discussion

It is however advisable to assert right away that the distribution
of a \va\ (random variable) $Q$ taking values in $\Q0$ must anyhow
be of a discrete type, allotting (possibly non vanishing)
probabilities to the individual rational numbers $q\in\Q0$:
conceivable continuous set functions -- namely with continuous,
albeit perhaps not absolutely continuous, \FD\ (cumulative
distribution function) -- would turn out to be not countably
additive, and hence would not qualify as measures, not to say as
probability distributions. Every continuous \FD\ for $Q$ would
indeed entail that at the same time $\pr{Q=q}=0,\;\forall q\in\Q0$,
and $\pr{Q\in\Q0}=1$, while $\Q0$ apparently is the countable union
of the disjoint, negligible sets $\{q\}$: in plain conflict with the
countable additivity. This in particular rules out for the numbers
in $\Q0$ also the possibility of being in some sense \emph{uniformly
distributed} (an imaginable surrogate of equiprobability suggested
by the rationals density): this property would in fact require for
$Q$ a \FD\ of the uniform type
\begin{equation*}
    F_Q(x)=\pr{Q\le x}=\left\{\begin{array}{ll}
                                       0 &\quad \hbox{$x<0$} \\
                                       x &\quad \hbox{$0\le x<1$} \\
                                       1 &\quad \hbox{$1\le x$}
                              \end{array}
                              \right.
\end{equation*}
which is apparently continuous, and would hence attribute
probability $0$ to every single $q$, but probability $1$ to $\Q0$

We would like to stress, moreover, that the problem focused on in
the present paper is not how to realistically produce -- possibly
equiprobable -- rational numbers \emph{at random}: this would be
performed in a trivial way, for instance, just by taking random,
uniformly distributed \emph{real} numbers, and then by truncating
them to a prefixed number $n$ of decimal digits, as always done in
practice in every computer simulations of random numbers in $[0,1]$.
It is apparent however that in so doing we would shrink $\Q0$ to a
\emph{finite} set of rational numbers (they would be exactly
$10^n+1$) that could always be made exactly equiprobable, failing
instead to allot a non vanishing probability to the remaining,
overwhelmingly more numerous, elements of $\Q0$. The aim of our
inquiry is instead to find a sensible way to attribute (non
vanishing, and possibly not too different from each other)
probabilities \emph{to every rational number} in $\Q0$, their
practical simulation being considered here but an eventual side
effect of this allocation

Remark that one could be lured to think that a way around the
previous snag could consist in drawing again uniformly distributed
real numbers, yet truncating the decimal digits to some
\emph{random} number $N$ taking arbitrary, finite but unbounded
integer values. Even in this way, however, not every rational number
would have a chance to be produced: the said procedure would indeed
\emph{a priori} exclude all the (infinitely many) rational numbers
with an infinite, periodic decimal representation, as for instance
$^1/_3,\,\!^2/_3,\ldots$ and so on. In the light of this preliminary
scrutiny the best way to tackle the task of laying down a
probability on $\Q0$ seems then to be to exploit the fractional
representation $q=\,\!^n/_m$ of every rational number by attributing
some suitable joint distribution to its numerators and denominators

\section{Distributions on $\Q0$}\label{ratiodistr}

Taking advantage of the well known diagram used to prove the
countability of the rational numbers, we will consider two dependent
\va's $M$ and $N$ with integer values
\begin{equation*}
    m=1,2,\ldots \qquad\qquad n=0,1,2,\ldots,m
\end{equation*}
and acting respectively as denominator and numerator of the random
rational number $Q=\,^N/_M\in[0,1]$. As a consequence $Q$ will take
the values $q=\,^n/_m$ arrayed in a triangular scheme as in
Table\myref{rationals}.
\begin{table}
\begin{center}
    \begin{tabular}{|c|cccccccccc}
  &  &  &  &  & $\bm n$ &  &  &  & & \\
  \hline
    $\bm m$ & $0$ & $1$ & $2$ & $3$ & $4$ & $5$ &$6$  & $7$ & $8$ & $\ldots$ \\
  \hline
  $1$  & $0$ & $1$ &  &  &  &  &  & & &\\
  $2$  & $0$ & $^1/_2$ & $1$ &  &  &  &  & & &\\
  $3$  & $0$ & $^1/_3$ & $^2/_3$ & $1$ &  &  &  & & &\\
  $4$  & $0$ & $^1/_4$ & $^2/_4$ & $^3/_4$ & $1$ &  &  & & &\\
  $5$  & $0$ & $^1/_5$ & $^2/_5$ & $^3/_5$ & $^4/_5$ & $1$ &  & & & \\
  $6$  & $0$ & $^1/_6$ & $^2/_6$ & $^3/_6$ & $^4/_6$ & $^5/_6$ & $1$ & & &\\
  $7$  & $0$ & $^1/_7$ & $^2/_7$ & $^3/_7$ & $^4/_7$ & $^5/_7$ & $^6/_7$ & $1$ & &\\
  $8$  & $0$ & $^1/_8$ & $^2/_8$ & $^3/_8$ & $^4/_8$ & $^5/_8$ & $^6/_8$ & $^7/_8$ & $1$ &\\
  $\vdots$  & $\vdots$ &  &  &  &  &  &  & & & $\ddots$\\
  %\hline
\end{tabular}
\caption{Table of rational numbers $q=\,\!^n/_m$ with repetitions:
many fractions are reducible to canonical forms already present in
earlier positions}\label{rationals}
\end{center}
\end{table}
It is apparent however that in this way every rational number $q$
shows up infinitely many times due to the presence of reducible
fractions: for instance -- with the usual notation for repeating
decimals -- we have
\begin{equation*}
    0.5=\,\!^1/_2=\,\!^2/_4=\,\!^3/_6=\ldots\qquad 0.\overline{3}=\,\!^1/_3=\,\!^2/_6=\ldots\qquad 0.75=\,\!^3/_4=\,\!^6/_8=\ldots
\end{equation*}
and hence, to avoid repetitions, the rational numbers in $[0,1]$
should rather be listed with blanks as in
Table\myref{norepetitions}. While always possible in principle,
however, it would be uneasy to assign probabilities directly to the
elements of the said Table\myref{norepetitions}: there is in fact no
simple way to attribute a progressive index to them (what for
instance is the $1\,000^{th}$ element?) since the numbers $\nu_m$ of
the different rationals in every row sharing a common irreducible
denominator $m$ constitute a rather irregular sequence, as we will
briefly discuss in the Section\myref{concl}. As a consequence it is
advisable to take advantage of the complete Table\myref{rationals}
by introducing a joint distributions of $N$ and $M$
\begin{table}
\begin{center}
    \begin{tabular}{c|c|ccccccccc}
   &  &  &  &  &  & $\bm n$ &  &  &  &  \\
  \hline
 $\bm{ \nu_m}$  & $\bm m$ & $0$ & $1$ & $2$ & $3$ & $4$ & $5$ &$6$  & $7$ & $\ldots$   \\
  \hline
2 &  $1$  & $0$ & $1$ &  &  &  &  &  & & \\
1 &  $2$  &  & $^1/_2$ &  &  &  &  &  & & \\
2 &  $3$  &  & $^1/_3$ & $^2/_3$ &  &  &  &  & & \\
2 &  $4$  &  & $^1/_4$ &  & $^3/_4$ &  &  &  & & \\
4 &  $5$  &  & $^1/_5$ & $^2/_5$ & $^3/_5$ & $^4/_5$ &  &  & &  \\
2 &  $6$  &  & $^1/_6$ &  &  &  & $^5/_6$ &  & & \\
6 &  $7$  &  & $^1/_7$ & $^2/_7$ & $^3/_7$ & $^4/_7$ & $^5/_7$ & $^6/_7$ &  & \\
4 &  $8$  &  & $^1/_8$ &  & $^3/_8$ &  & $^5/_8$ &  & $^7/_8$ &  \\
$\vdots$ &  $\vdots$  &  & $\vdots$ &  &  &  &  &  & & $\ddots$  \\
  %\hline
\end{tabular}
\caption{Table of rational numbers $q\doteq\,\!^n/_m$ without
repetitions: only irreducible fractions are represented, along with
the number $\nu_m$ of the different rationals sharing a common
irreducible denominator $m$}\label{norepetitions}
\end{center}
\end{table}
\begin{align*}
  &  \pr{M=m} \qquad\qquad\qquad\; m=1,2,\ldots \\
  &  \pr{N=n\left|M=m\right.} \qquad\quad \,n=0,1,2,\ldots,m \\
  &  \pr{N=n,M=m}=\pr{N=n\left|M=m\right.}\,\pr{M=m}
\end{align*}
For a rational number $q$ we will also adopt the notation
\begin{equation*}
    q\doteq\,\!^n/_m
\end{equation*}
to indicate that $^n/_m$ is the irreducible representation of $q$,
namely that $n$ and $m$ are co-primes: for instance in the previous
examples it will be
\begin{equation*}
    0.5\doteq\,\!^1/_2\qquad\quad0.\overline{3}\doteq\,\!^1/_3\qquad\quad0.75\doteq\,\!^3/_4
\end{equation*}
For every rational $q\doteq\,\!^n/_m$ we will then have the discrete
distribution
\begin{eqnarray}\label{discrD}
    \pr{Q=q}&=&\sum_{\ell=1}^\infty\pr{N=\ell n,M=\ell m}\nonumber\\
  &=&\sum_{\ell=1}^\infty\pr{N=\ell n\left|M=\ell m\right.}\,\pr{M=\ell m}
\end{eqnarray}
which gives a probability to every rational number $0\le q\le1$.
This also allows to define the \FD\ of $Q$ as (here of course
$x\in\RE$)
\begin{eqnarray}
  F_Q(x) &=& \pr{Q\le x}=\pr{N\le Mx}=\sum_{m=1}^\infty\pr{N\le
  mx\left|M=m\right.}\pr{M=m}\nonumber\\
   &=&\sum_{m=1}^\infty F_N(mx|M=m)\pr{M=m}\label{cdf1}
\end{eqnarray}
and hence also the probability of $Q$ falling in $(a,b]$ for $0\le
a<b\le1$ real numbers:
\begin{eqnarray}
  \pr{a< Q\le b}&=&F_Q(b)-F_Q(a)\nonumber\\
  &=&\sum_{m=1}^\infty \big[F_N(mb|M=m)-F_N(ma|M=m)\big]\pr{M=m}\label{probint1}
\end{eqnarray}
Notice that the conditional \FD\ of $N$ can also be given as
\begin{eqnarray}
  F_N(x|M=m)&=&\pr{N\le
  x\left|M=m\right.}=\sum_{n=0}^m\pr{N=n\left|M=m\right.}\vartheta(x-n)\nonumber\\
   &=&\sum_{n=0}^{\lfloor
   x\rfloor}\pr{N=n\left|M=m\right.}\label{condcdf}
\end{eqnarray}
where
\begin{equation*}
    \vartheta(x)=\left\{
                   \begin{array}{ll}
                     1 & \qquad x\ge0 \\
                     0 & \qquad x<0
                   \end{array}
                 \right.
\end{equation*}
is the Heaviside function, while for every real number $x$, the
symbol $\lfloor x\rfloor$ denotes the \emph{floor} of $x$, namely
the greatest integer less than or equal to $x$. As a consequence the
equations\refeq{cdf1} and\refeq{probint1} also take the
form\begin{eqnarray}
  F_Q(x)\!\! &=&\!\! \sum_{m=1}^\infty \pr{M=m}\sum_{n=0}^{\lfloor
   mx\rfloor}\pr{N=n\left|M=m\right.} \label{cdf}\\
  \pr{a< Q\le b}\!\! &=&\!\! \sum_{m=1}^\infty\pr{M=m}(1-\delta_{\lfloor ma\rfloor,\lfloor mb\rfloor})\!\!\! \sum_{n=\lfloor
   ma\rfloor+1}^{\lfloor
   mb\rfloor}\!\!\!\!\pr{N=n\left|M=m\right.} \label{probint}
\end{eqnarray}
where the Kronecker delta takes into account the fact that when
$\lfloor mb\rfloor=\lfloor ma\rfloor$ the term vanishes, so that
$\lfloor mb\rfloor\ge\lfloor ma\rfloor+1$. We moreover have for the
expectations and the characteristic function
\begin{eqnarray}
  \att{Q} &=& \att{\,\!^N/_M} =\att{\frac{1}{M}\att{N|M}}=\sum_{m=1}^\infty\frac{\pr{M=m}}{m}\att{N|M=m}\label{exp}\\
   \att{Q^2}&=&\sum_{m=1}^\infty\frac{\pr{M=m}}{m^2}\att{\left.N^2\right|M=m}\label{exp2} \\
   \varphi_Q(u)&=&\att{e^{iuN/M}}=\sum_{m=1}^\infty\pr{M=m}\sum_{n=0}^m
   e^{iun/m}\pr{N=n\left|M=m\right.}\nonumber \\
   &=&\sum_{m=1}^\infty\pr{M=m}\varphi_N\left(^u/_m\left|M=m\right.\right)\label{chf}
\end{eqnarray}
where we also adopted the shorthand notation
\begin{equation*}
    \varphi_N\left(u\left|M=m\right.\right)=\att{e^{iuN}\left|M=m\right.}=\sum_{n=0}^m
   e^{iun}\pr{N=n\left|M=m\right.}
\end{equation*}
The actual joint distributions of $N$ and $M$ can now be chosen in
several ways, and we go on now in the next sections to survey a few
particular cases

\section{Equiprobable numerators}\label{equinum}

Let us suppose now for simplicity that for a given denominator
$m\ge1$ the $m+1$ possible values of the numerator $n=0,1,\ldots,m$
are equiprobable in the sense that
\begin{equation*}
  \pr{N=n\left|M=m\right.}=\frac{1}{m+1}\qquad\quad n=0,1,\ldots,m
\end{equation*}
We then have (see\mycite{grad} $\bm{0.121}$)
\begin{eqnarray*}
  \att{N|M=m} &=& \sum_{n=0}^m\frac{n}{m+1}=\frac{1}{m+1}\,\frac{m(m+1)}{2}=\frac{m}{2} \\
  \att{\left.N^2\right|M=m} &=&\sum_{n=0}^m\frac{n^2}{m+1}=\frac{1}{m+1}\,\frac{m(m+1)(2m+1)}{6}=\frac{m(2m+1)}{6}
\end{eqnarray*}
and hence from\refeq{exp} and\refeq{exp2}
\begin{eqnarray*}
  \att{Q} &=& \sum_{m=1}^\infty\frac{\pr{M=m}}{m}\,\frac{m}{2}=\frac{1}{2}\sum_{m=1}^\infty\pr{M=m} =\frac{1}{2}\\
  \att{Q^2} &=&\sum_{m=1}^\infty\frac{\pr{M=m}}{m^2}\att{N^2|M=m}=\sum_{m=1}^\infty\frac{2m+1}{6m}\pr{M=m}\\
  &=&\frac{1}{3}\sum_{m=1}^\infty\pr{M=m}+\frac{1}{6}\sum_{m=1}^\infty\frac{\pr{M=m}}{m}=\frac{1}{3}+\frac{1}{6}\att{\frac{1}{M}}\\
  \var{Q}&=&\att{Q^2}-\att{Q}^2=\frac{1}{12}+\frac{1}{6}\att{\frac{1}{M}}
\end{eqnarray*}
As for the distribution, with $n,m$ co-primes and $0\le n\le m$,
from\refeq{discrD} we have
\begin{equation}\label{ratprob}
  \pr{Q=q} =
  \sum_{\ell=1}^\infty\frac{\pr{M=\ell m}}{\ell m+1}\qquad\quad q\doteq\,\!^n/_m
\end{equation}
which is apparently independent from $n$ and is contingent only on
the value of the irreducible denominator $m$. The characteristic
function\refeq{chf} moreover is
\begin{equation*}
  \varphi_N\left(u\left|M=m\right.\right)= \frac{1}{m+1}\sum_{n=0}^me^{iun}
  \qquad\quad
  \varphi_Q(u) =\sum_{m=1}^\infty\frac{\pr{M=m}}{m+1}\sum_{n=0}^me^{iun/m}
\end{equation*}
while for the \FD\refeq{cdf} we have from\refeq{condcdf}
\begin{eqnarray}
  F_N\left(mx\left|M=m\right.\right)&=&
  \frac{1}{m+1}\sum_{n=0}^m\vartheta\left(mx-n\right)=\left\{
                                                                 \begin{array}{ll}
                                                                   0 &\quad \hbox{$x<0$} \\
                                                                   \frac{\lfloor mx\rfloor+1}{m+1} &\quad \hbox{$0\le x<1$} \\
                                                                   1 &\quad \hbox{$1\le x$}
                                                                 \end{array}
                                                               \right.\nonumber
  \\
  F_Q(x) &=&\sum_{m=1}^\infty\frac{\pr{M=m}}{m+1}\sum_{n=0}^m\vartheta\left(mx-n\right)\nonumber\\
                                                               &=&\left\{
                                                                 \begin{array}{ll}
                                                                   0 &\quad \hbox{$x<0$} \\
                                                                   \sum_{m\ge1}\pr{M=m}\frac{\lfloor mx\rfloor+1}{m+1} &\quad \hbox{$0\le x<1$} \\
                                                                   1 &\quad \hbox{$1\le x$}
                                                                 \end{array}
                                                               \right.\label{ratcdf}
\end{eqnarray}
and the probability\refeq{probint} with $0\le a<b\le1$ becomes
\begin{equation}\label{ratint}
  \pr{a< Q\le b}=\sum_{m=1}^\infty\pr{M=m}\frac{\lfloor mb\rfloor-\lfloor ma\rfloor}{m+1}
\end{equation}
It is apparent then that -- but for the value of the expectation
$\att{Q}=\,\!^1/_2$ -- all these quantities depend on the choice of
the denominator distribution. That notwithstanding we will show in
the next section that, under reasonable conditions on the
denominators $M$, the distribution of $Q$ can in fact be made as
near as we want to -- but not exactly coincident with -- a uniform
distribution in $[0,1]$: a behavior that we dubbed \emph{asymptotic
equiprobability}

\section{Asymptotic equiprobability}\label{equi}

For the equiprobable numerators introduced in the previous section,
and by denoting for short as $p_m=\pr{M=m}$ the distribution of $M$,
and as $s=\sup_m p_m$ the supremum of all its values, let us take
now a sequence of denominators $\{M_k\}_{k\ge1}$ with distributions
$\{p_m(k)\}_{k\ge1}$, and with $s_k$ vanishing for $k\to\infty$ in
such a way that
\begin{equation}\label{log}
    \lim_k\, s_k\ln k=0
\end{equation}
In other words we consider a sequence of distributions that are
increasingly (and uniformly) flattened toward zero, so that the
denominators too are increasingly equiprobable. Ready examples of
these sequences with $k=1,2,\ldots$ are for instance that of the
\emph{finite equiprobable} distributions
\begin{equation*}
    p_m(k)=\left\{
             \begin{array}{ll}
               ^1/_k & \hbox{$m=1,2,\ldots,k$} \\
               0 & \hbox{$m>k$}
             \end{array}
           \right.
\end{equation*}
where apparently $s_k=\,\!^1/_k\tok0\,$; that of the
\emph{geometric} distributions
\begin{equation*}
    p_m(k)=w_k(1-w_k)^{m-1}\quad\qquad m=1,2,\ldots
\end{equation*}
with infinitesimal $w_k$ so that $s_k=w_k\tok0\,$; and finally that
of the \emph{Poisson} distributions
\begin{equation*}
    p_m(k)=e^{-\lambda_k}\frac{\lambda_k^{m-1}}{(m-1)!}\quad\qquad m=1,2,\ldots
\end{equation*}
with divergent $\lambda_k$, where again the modal values are
infinitesimal: we know indeed that a Poisson distribution
%\begin{equation*}
%    e^{-\lambda}\frac{\lambda^{m-1}}{(m-1)!}=e^{-\lambda}\frac{\lambda^{m-1}}{\Gamma(m)}\quad\qquad m=1,2,\ldots
%\end{equation*}
attains its maximum in $\lfloor\lambda_k\rfloor+1$, so that for
$\lambda_k\tok+\infty$ its modal value $s_k$ essentially behaves as
(see\mycite{grad} $\bm{8.327.1}$)
\begin{equation*}
    s_k=e^{-\lambda_k}\frac{\lambda_k^{\lambda_k-1}}{\Gamma(\lambda_k)}=\frac{1}{\sqrt{2\pi\lambda_k}\big(1+{O}(\lambda_k^{-1})\big)}\tok0
\end{equation*}
 \begin{lem}
Within the previous notations and conditions we have
\begin{equation}\label{lem}
   \mu_k=\att{\,\!^1/_{M_k}}=\sum_{m=1}^{\infty}\frac{p_m(k)}{m}\tok0
\end{equation}
 \end{lem}
 \indim
The positive series defining $\mu_k$ is certainly convergent because
\begin{equation*}
    \mu_k=\sum_{m=1}^{\infty}\frac{p_m(k)}{m}<\sum_{m=1}^{\infty}p_m(k)=1
\end{equation*}
and hence we can always write
\begin{equation*}
    \mu_k=\sum_{m=1}^{\infty}\frac{p_m(k)}{m}=\sum_{m=1}^{k}\frac{p_m(k)}{m}+R_k
\end{equation*}
where
\begin{equation*}
    R_k=\sum_{m=k+1}^{\infty}\frac{p_m(k)}{m}\tok0
\end{equation*}
is an infinitesimal remainder. Remark that now $k$ plays both the
roles of index of the distribution sequence, and of series cut-off.
On the other hand, under our stated conditions
\begin{equation*}
    \sum_{m=1}^{k}\frac{p_m(k)}{m}<s_k\sum_{m=1}^{k}\frac{1}{m}=s_k
    H_k
\end{equation*}
where $H_k$ denotes the $k^{th}$ \emph{harmonic number}, namely the
sum of the reciprocal integers up to $^1/_k$: it is well\ known
(\cite{grad} $\bm{0.131}$) that for $k\to\infty$ the $H_k$ grow as
$\ln k$, so that from\refeq{log} we have $s_kH_k\tok0$, and finally
$\mu_k=s_kH_k+R_k\tok0$
 \findim
 \begin{prop}
If $Q=\,\!^N/_M$ and $F_Q(x)$ is its \emph{\FD}, then, within the
notation and conditions outlined above, we have
\begin{align}
  & \lim_{\mmax}\pr{Q=q}=0\qquad\qquad\lim_{\mmax}\pr{a<Q\le
               b}=b-a\label{equiprob}\\
%\end{equation}
%and in particular
%\begin{equation}
   & \qquad\qquad\quad\lim_{\mmax}F_Q(x)=\left\{\begin{array}{ll}
                                       0 &\quad \hbox{$x<0$} \\
                                       x &\quad \hbox{$0\le x<1$} \\
                                       1 &\quad \hbox{$1\le x$}
                              \end{array}
                              \right.\label{unifcdf}
\end{align}
 \end{prop}
 \indim
Since our series have positive terms the first result
in\refeq{equiprob} follows from\refeq{ratprob} and\refeq{lem}
because, with $q\doteq\,\!^n/_j$
\begin{equation*}
    \pr{Q=q}=\sum_{\ell=1}^\infty\frac{p_{\ell j}(k)}{\ell j+1}<\sum_{m=1}^\infty\frac{p_{m}(k)}{m+1}
    <\sum_{m=1}^\infty\frac{p_{m}(k)}{m}=\mu_k\tok0
\end{equation*}
As for the second result in\refeq{equiprob}, since for every real
number $x$ it is $x-1\le\lfloor x\rfloor\le x$, for every
$\mmax=1,2,\ldots$, and $0\le a<b\le1$, we have from\refeq{ratint}
\begin{equation*}
   \sum_{m=1}^{\infty}p_m(k)\frac{m(b-a)-1}{m+1} \le\,\pr{a<Q\le b}\,\le\sum_{m=1}^{\infty}p_m(k)\frac{m(b-a)+1}{m+1}
\end{equation*}
namely
\begin{equation*}
   b-a+(a-b-1)\sum_{m=1}^{\infty}\frac{p_m(k)}{m+1}\le\,\pr{a<Q\le b}\,\le b-a+(a-b+1)\sum_{m=1}^{\infty}\frac{p_m(k)}{m+1}
\end{equation*}
so that, since $a-b-1\le0$ and $a-b+1\ge0$, it is
\begin{equation*}
    b-a+(a-b-1)\mu_k\le\,\pr{a<Q\le b}\,\le b-a+(a-b+1)\mu_k
\end{equation*}
The second result\refeq{equiprob} follows then from\refeq{lem}. In a
similar way we finally find for\refeq{unifcdf} that
\begin{equation*}
   \sum_{m=1}^{\infty}p_m(k)\frac{mx}{m+1} \le F_Q(x)\le\sum_{m=1}^{\infty}p_m(k)\frac{mx+1}{m+1}\qquad\quad 0\le x\le1
\end{equation*}
namely
\begin{equation*}
   x-x\sum_{m=1}^{\infty}\frac{p_m(k)}{m+1}\le F_Q(x)\le x+(1-x)\sum_{m=1}^{\infty}\frac{p_m(k)}{m+1}
\end{equation*}
and hence
\begin{equation*}
   x-x\,\mu_k< F_Q(x)< x+(1-x)\mu_k
\end{equation*}
so that the result again follows from\refeq{lem}
 \findim

\noindent From this proposition we see that in the limit
$\mmax\to\infty$, while the probability of every single rational
number rightly vanishes, the probability of these numbers lumped
together in intervals does not: a behavior highly reminiscent of
what happens to continuously distributed \emph{real} \va's. For the
reasons presented in the Section\myref{ratio}, however, the previous
result by no means imply that we can implement a uniform limit
distribution on $\Q0$ (as we said, there is not such a thing), but
it rather suggests that our random rational numbers $Q$ -- at least
for denominators $m$ distributed in a fairly flat way, and
numerators $n$ conditionally equiprobable between $0$ and $m$ --
asymptotically behave as uniformly distributed in $[0,1]$, and hence
they quite reasonably correspond to our intuitive idea of
\emph{taking rational numbers at random.} In this perspective remark
also that, under our conditions, we have for the variance
\begin{equation*}
    \var{Q}=\frac{1}{12}+\frac{1}{6}\,\att{\frac{1}{M}}=\frac{1}{12}+\frac{\mu_k}{6}\tok\frac{1}{12}
\end{equation*}
again in agreement with an approximate uniform distribution in
$[0,1]$

\section{Denominator distributions}\label{denom}

\subsection{Geometric denominators}\label{geometric}
 \begin{figure}
 \begin{center}
\includegraphics*[width=15cm]{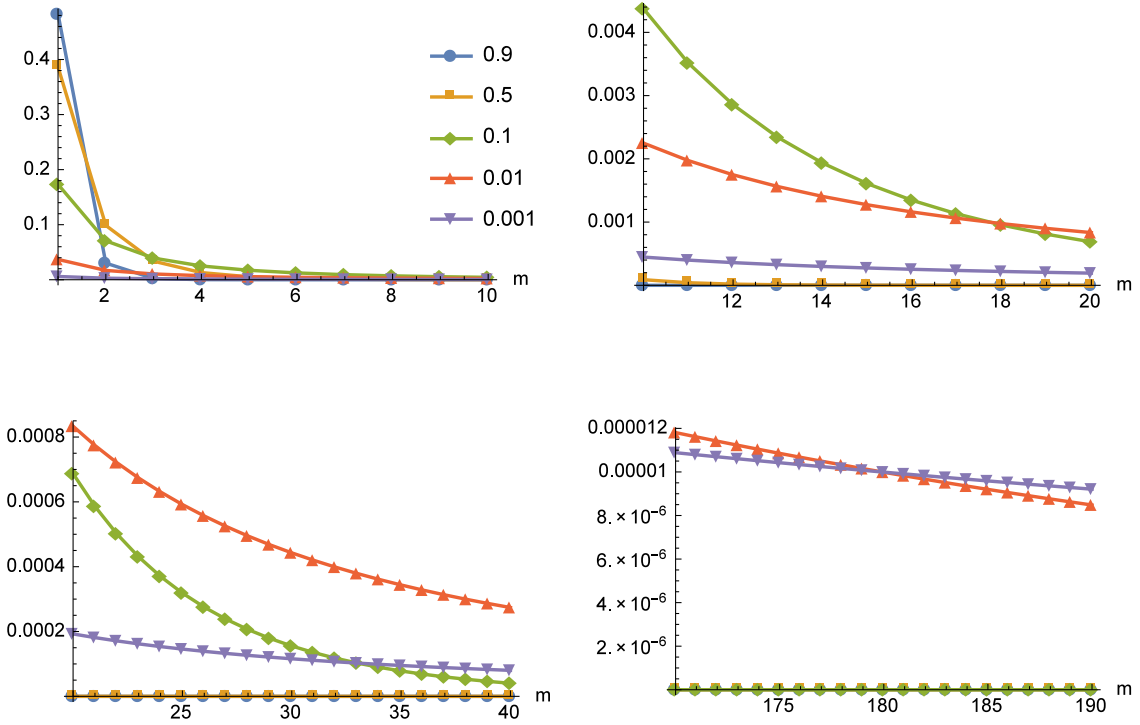}
\caption{Probabilities attributed to rational numbers as a function
of the irreducible, geometrically distributed denominators $m$, and
for decreasing ($0.9, 0.5, 0.1, 0.01, 0.001$) values of $w$: by
choosing different $m$ intervals, the pictures show how these
probabilities level down to infinitesimal equiprobability for
$w\to0$ }\label{EquiprobGeometric}
 \end{center}
 \end{figure}

A few closed formulas about the \va\ $Q$ are available for
particular denominator distributions: let us suppose for instance
that the denominator $M$ be geometrically distributed as
\begin{equation*}
    \pr{M=m}=w\,(1-w)^{m-1}\qquad\quad w>0,\qquad m=1,2,\ldots
\end{equation*}
In this case we first find
\begin{equation*}
    \att{M}=\frac{1}{w}\quad\qquad\att{\frac{1}{M}}=\sum_{m=1}^\infty\frac{w(1-w)^{m-1}}{m}=\frac{w}{1-w}\sum_{m=1}^\infty\frac{(1-w)^m}{m}=-\frac{w\,\ln w}{1-w}
\end{equation*}
and hence
\begin{equation*}
  \var{Q}=\frac{1}{12}-\frac{w\,\ln w}{6(1-w)}
\end{equation*}
while for the \FD\ we do not go beyond its formal definition
\begin{equation*}
    F_Q(x) =\left\{ \begin{array}{ll}                    0 &\quad \hbox{$x<0$} \\
                                                                   \sum_{m\ge1}w\,(1-w)^{m-1}\frac{\lfloor mx\rfloor+1}{m+1} &\quad \hbox{$0\le x<1$} \\
                                                                   1 &\quad \hbox{$1\le x$}
                                                                 \end{array}
                                                               \right.
\end{equation*}
As for the $Q$ distribution instead, taking $q\doteq\,\!^n/_m$, we
find with $j=\ell-1$ the analytic expression
\begin{eqnarray}
 \pr{Q=q}  &=& \sum_{\ell=1}^\infty\frac{w(1-w)^{\ell m-1}}{\ell m+1}=w(1-w)^{m-1}\sum_{\ell=1}^\infty\frac{(1-w)^{m(\ell-1)}}{\ell m+1}\nonumber \\
   &=&w(1-w)^{m-1}\sum_{j=0}^\infty\frac{(1-w)^{mj}}{m(j+1)+1}\nonumber\\
   &=&\frac{w(1-w)^{m-1}}{m+1}\,_2F_1\left(1,1+\,\!^1/_m\,;2+\,\!^1/_m\,;(1-w)^m\right)\label{geomdistr}
\end{eqnarray}
where $_2F_1(a,b;c;z)$ is a hypergeometric function~\cite{grad} that
gauges the deviation of $\pr{Q=q}$ from the corresponding joint
probability of $N,M$
\begin{equation*}
    \pr{N=n,M=m}=\frac{w(1-w)^{m-1}}{m+1}
\end{equation*}
This formula allows a graphic representation of $\pr{Q=q}$ as a
function of the irreducible denominators $m$ displayed in the
Figure\myref{EquiprobGeometric} where it is apparent how the initial
($m=1$) ordering of the probabilities (increasing with the $w$
values going from $w=0.001$ to $w=0.9$) becomes totally overturned
for $m$ large enough. Remark that each value of the
probability\refeq{geomdistr} should be understood as attributed to
every rational number with the same $m$ as irreducible denominator,
for instance (see Table\myref{norepetitions}): for $m=1$ we get the
probability of $q\doteq0$ and $1$; for $m=2$ the probability of
$q\doteq\,\!^1/_2$ alone; for $m=3$ that of
$q\doteq\,\!^1/_3\,,^2/_3\,$; for $m=4$ that of
$q\doteq\,\!^1/_4\,,^3/_4\,;\,\ldots\,$ and so on. This allows, in
particular, to steer clear of an easy misunderstanding: it must be
noticed indeed that, while apparently
\begin{equation*}
    \sum_{m=1}^\infty\sum_{n=0}^m\pr{N=n,M=m}=\sum_{m=1}^\infty\sum_{n=0}^m\frac{w(1-w)^{m-1}}{m+1}=\sum_{m=1}^\infty
w(1-w)^{m-1}=1
\end{equation*}
we find instead
\begin{equation*}
    \sum_{m=1}^\infty\sum_{n=0}^m\frac{w(1-w)^{m-1}}{m+1}\,_2F_1\left(1,1+\,\!^1/_m\,;2+\,\!^1/_m\,;(1-w)^m\right)<1
\end{equation*}
as can be seen from the fact that for $0<w<1$
\begin{equation*}
   _2F_1\left(1,1+\,\!^1/_m\,;2+\,\!^1/_m\,;(1-w)^m\right)\;\left\{
                                                             \begin{array}{ll}
                                                               =1 & \qquad\hbox{$m=1$} \\
                                                               <1 & \qquad\hbox{$m=2,3,\ldots$}
                                                             \end{array}
                                                           \right.
\end{equation*}
This however is not in contradiction with the mandatory requirement
that
 \begin{figure}
 \begin{center}
\includegraphics*[width=12cm]{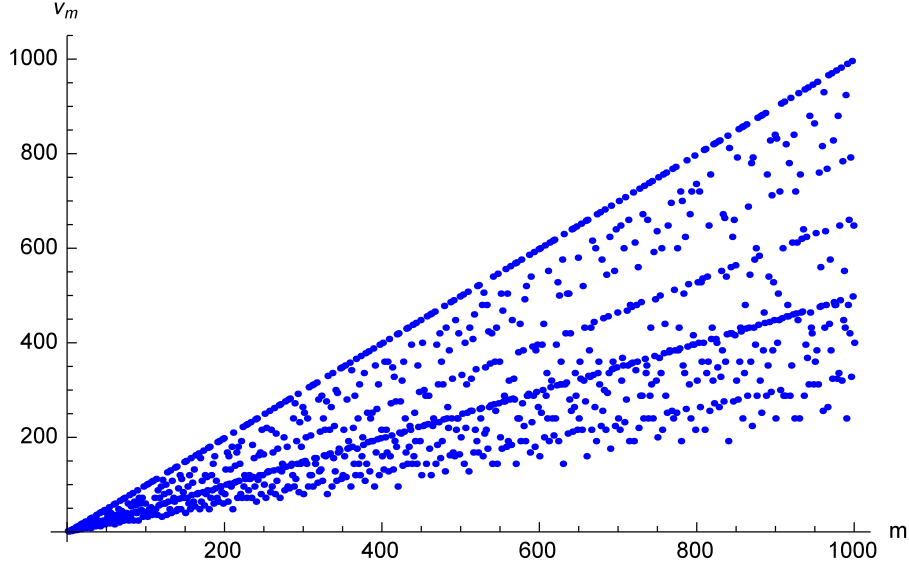}
\caption{Numerosity $\nu_m$ of the different rational numbers
$q\doteq\,\!^n/_m$ sharing a common, irreducible denominator
$m$}\label{DifferentRat}
 \end{center}
 \end{figure}
\begin{equation}\label{norm}
    \sum_{q\in\Q0}\pr{Q=q}=1
\end{equation}
precisely because -- as previously remarked -- the probability
associated to an $m$ must be attributed to several different
rational numbers $q$: if $\nu_m$ is the number of rationals $q$ that
have $m$ as its irreducible denominator, then we should rather pay
attention to ascertain the normalization in the form
\begin{eqnarray*}
   \sum_{q\in\Q0}\pr{Q=q}&=& \sum_{m=1}^\infty\sum_{n=0}^m\nu_m\,\frac{w(1-w)^{m-1}}{m+1}\,_2F_1\left(1,1+\,\!^1/_m\,;2+\,\!^1/_m\,;(1-w)^m\right)\\
                                  &=& \sum_{m=1}^\infty\nu_m\,w(1-w)^{m-1}\,_2F_1\left(1,1+\,\!^1/_m\,;2+\,\!^1/_m\,;(1-w)^m\right) =1
\end{eqnarray*}
Yet this result -- that we can consider as secured by construction
and definition -- is not easy to check by direct calculation because
a closed form for the sequence $\nu_m$ is not readily available: its
behavior is indeed rather irregular, albeit on average steadily
growing, as can be seen from an empirical plot of its first values
displayed in the Figure\myref{DifferentRat}. We postpone to the
Section\myref{concl} a few additional remarks about this point
showing in particular how the previous normalization condition can
instead be used to sequentially calculate the values of $\nu_m$

\subsection{Poisson and equiprobable denominators}\label{PoissEqui}

When on the other hand the denominators are distributed according to
other (albeit simple) laws we unfortunately no longer find
elementary closed forms for $\pr{Q=q}$. If the for instance $M$ is
Poisson distributed as
\begin{equation*}
    \pr{M=m}=e^{-\lambda}\,\frac{\lambda^{m-1}}{(m-1)!}\qquad\quad \lambda>0,\qquad m=1,2,\ldots
\end{equation*}
we find $\att{M}=1+\lambda$ and
\begin{equation*}
    \att{\frac{1}{M}}=\sum_{m=1}^\infty\frac{e^{-\lambda}}{m}\,\frac{\lambda^{m-1}}{(m-1)!}=\frac{e^{-\lambda}}{\lambda}\sum_{m=1}^\infty\frac{\lambda^m}{m!}
    =\frac{e^{-\lambda}}{\lambda}\left(\sum_{m=0}^\infty\frac{\lambda^m}{m!}-1\right)=\frac{1-e^{-\lambda}}{\lambda}
\end{equation*}
while for the variance we have
\begin{equation*}
  \var{Q}=\frac{1}{12}+\frac{1-e^{-\lambda}}{6\lambda}
\end{equation*}
but the \FD\ is
\begin{equation*}
    F_Q(x) =\left\{ \begin{array}{ll}
                                        0 &\quad \hbox{$x<0$} \\
                                    e^{-\lambda}\,\sum_{m\ge1}\frac{\lambda^{m-1}}{(m-1)!}\,\frac{\lfloor mx\rfloor+1}{m+1} &\quad \hbox{$0\le x<1$} \\
                                        1 &\quad \hbox{$1\le x$}
                                   \end{array}
                                   \right.
\end{equation*}
and for the distribution, taking $q\doteq\,\!^n/_m$, we find
\begin{equation*}
 \pr{Q=q}  = \sum_{\ell=1}^\infty\frac{e^{-\lambda}}{\ell m+1}\,\frac{\lambda^{\ell m-1}}{(\ell m-1)!}
\end{equation*}
with no closed expression readily available
 \begin{figure}
 \begin{center}
\includegraphics*[width=12cm]{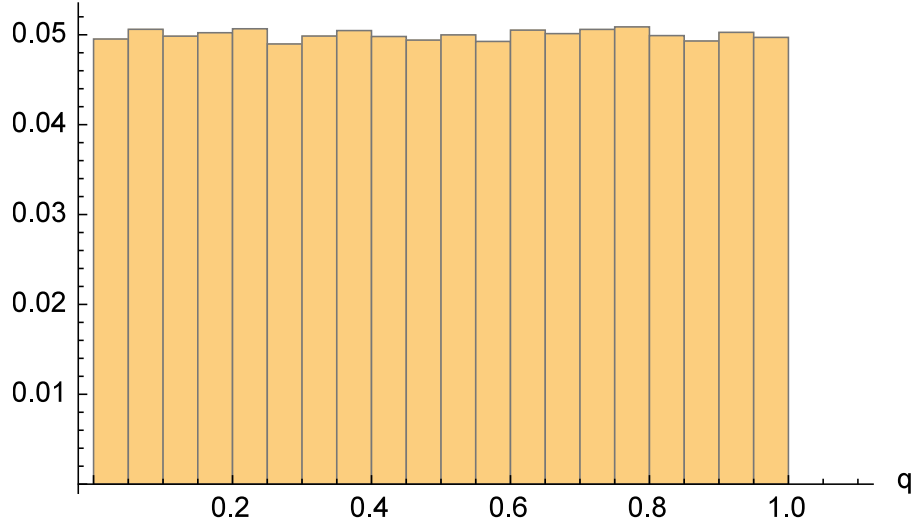}
\caption{Typical histogram of the relative frequencies of a sample
of $10^5$ random rationals generated following the procedure
described in the Section\myref{PoissEqui}: here the maximum value of
the equiprobable denominators is chosen to be $k=10^5$
}\label{EquiHisto}
 \end{center}
 \end{figure}

Consider instead denominators $M$ taking only a finite number
$\mmax=1,2,\ldots$ of equiprobable values $m$:
\begin{equation*}
    \pr{M=m}=\left\{
             \begin{array}{ll}
               ^1/_k & \hbox{$m=1,2,\ldots,k$} \\
               0 & \hbox{$m>k$}
             \end{array}
           \right.
\end{equation*}
We then have
\begin{equation*}
    \att{M}=\frac{k+1}{2}\qquad\qquad
    \att{\frac{1}{M}}=\frac{1}{\mmax}\sum_{m=1}^{\mmax}\frac{1}{m}=\frac{H_{\mmax}}{\mmax}
\end{equation*}
and hence
\begin{equation*}
  \var{Q} = \frac{1}{12}+\frac{H_{\mmax}}{6\mmax}
\end{equation*}
while for the \FD\ it is
\begin{equation*}
    F_Q(x)=\frac{1}{\mmax}\sum_{m=1}^{\mmax}\frac{\lfloor mx\rfloor+1}{m+1}
\end{equation*}
and the discrete distribution probabilities are
\begin{equation*}
    \pr{Q=q}=\frac{1}{k}\sum_{\ell=1}^{\lfloor\,\!^k/_m\rfloor}\frac{1}{\ell
m+1}\qquad\quad q\doteq\,\!^n/_m
\end{equation*}
where, since $m\le k$, it is always
$\lfloor\,\!^k/_m\rfloor=1,2,\ldots$.

Even in this case we have then no closed formulas to show, but since
the sums involved are now always finite this seems to hint to a
simple -- but essentially trivial -- procedure to simulate an
\emph{asymptotically equiprobable sample} of rational numbers in
$[0,1]$: choose first a large enough number $\mmax$, then sample a
random integer $m$ among the equiprobable numbers
$1,2,\ldots,\mmax$, and finally a random integer $n$ among the
equiprobable numbers $0,1,\ldots,m$ and put $q=\,\!^n/_m$. By
repeating this procedure a number of times large enough we get a
sample almost uniformly distributed in $[0,1]$ as shown in
Figure\myref{EquiHisto}. The drawback of this procedure, however, as
already remarked in the Section\myref{ratio}, is that not every
number in $\Q0$ would have the chance of being drawn because only a
finite number among them would actually be taken into account. This
large but finite set of numbers could also be made, in principle,
exactly equiprobable, but the infinitely many remaining rational
numbers would instead be totally excluded with strictly zero
probability

\section{Sequencing rational numbers}\label{concl}

Other examples of distributions on the rational numbers in $[0,1]$
are of course possible: for instance, with $0<p<1$ and for a given
denominator $m=1,2,\ldots$, it is possible to suppose that the
numerators are binomially -- instead of equiprobably -- distributed
as
\begin{equation*}
    \pr{N=n\left|M=m\right.}={m\choose n}p_ì(1-p)^{m-n}\qquad\quad n=0,1,\ldots,m
\end{equation*}
By choosing then a suitable distribution for the denominator $M$ we
can define the global distribution of $Q=\,\!^N/_M$. However, rather
than indulging in displaying these further examples, we would like
to conclude this paper with a few remarks about a particular
residual open problem

%by remarking that, even when the $Q$ distribution is very similar to
%a continuous, uniform distribution on $[0,1]$, it remains a
%\emph{discrete distribution}

%\begin{appendix}

%\section{Sequencing rational numbers}\label{sequence}

\begin{table}
\begin{center}
    \begin{tabular}{c|c|c|ccccccccccc}
&  &  &  &  &  &  &  & $\bm n$ &  &  & &  &     \\
  \hline
$\bm{ \sigma_m}$ & $\bm{ \nu_m}$  & $\bm m$ & $0$ & $1$ & $2$ & $3$ & $4$ & $5$ &$6$  & $7$ & $8$ & $9$ & $10$    \\
  \hline
$1$ & $2$ &  $1$  & $0$ & $1$ &  &  &  &  &  & & &  &    \\
$^1/_2$ & $1$ &  $2$  &  & $^1/_2$ &  &  &  &  &  & & &  &    \\
$1$ & $2$ &  $3$  &  & $^1/_3$ & $^2/_3$ &  &  &  &  & & &  &    \\
$1$ & $2$ &  $4$  &  & $^1/_4$ &  & $^3/_4$ &  &  &  &  \\
$2$ & $4$ &  $5$  &  & $^1/_5$ & $^2/_5$ & $^3/_5$ & $^4/_5$ &  &  & &  &  &    \\
$1$ & $2$ &  $6$  &  & $^1/_6$ &  &  &  & $^5/_6$ &  & & &  &    \\
$3$ & $6$ &  $7$  &  & $^1/_7$ & $^2/_7$ & $^3/_7$ & $^4/_7$ & $^5/_7$ & $^6/_7$ &  & &  &    \\
$2$ & $4$ &  $8$  &  & $^1/_8$ &  & $^3/_8$ &  & $^5/_8$ &  & $^7/_8$ & &  &     \\
$3$ & $6$ &  $9$  &  & $^1/_9$ &  $^2/_9$ &  & $^4/_9$ & $^5/_9$ &  & $^7/_9$ & $^8/_9$ &  &       \\
$2$ & $4$ &  $10$  &  & $^1/_{10}$ &  & $^3/_{10}$ &  &  &  & $^7/_{10}$ & & $^9/_{10}$ &       \\
$5$ & $10$ &  $11$  &  & $^1/_{11}$ & $^2/_{11}$ & $^3/_{11}$ & $^4/_{11}$ & $^5/_{11}$ & $^6/_{11}$ & $^7/_{11}$ &$^8/_{11}$ &$^9/_{11}$  &$^{10}/_{11}$       \\
$\vdots$ & $\vdots$ &  $\vdots$  &  & $\vdots$ &  &  &  &  &  & &  &  &  $\vdots$  \\
  %\hline
\end{tabular}
\caption{Table of rational numbers $q\doteq\,\!^n/_m$ without
repetitions, along with the progressive number $\nu_m$ of the
different rationals sharing a common irreducible denominator $m$,
and their sums $\sigma_m$}\label{apptable}
\end{center}
\end{table}

We said from the beginning that since $\Q0$ is countable its
elements $0\le q\le1$ can certainly be arranged in a sequence. If on
the other hand we can manage to have in this sequence all the
rational numbers without repetitions, this would greatly facilitate
the task of giving a distribution on $\Q0$. In order however to put
in a sequence without repetitions $q_k$ all these rational numbers
in $[0,1]$ -- as listed for instance in the triangular, infinite
Table\myref{apptable} -- we should at least be able to find
regularities in their arrangement allowing to say with relative easy
both what is the rational number $q$ associated to an arbitrary
given index $k$, and viceversa what is the place (index $k$) of an
arbitrary given rational number $q$. But this quest is baffled by
the rather irregular running of the entries in the said triangular
table where, for instance, even the occurrence among the
denominators $m$ of the prime numbers (the only ones identifying
rows with no blanks beyond the extremes) is famously not immediately
predictable. It is apparent however that the possibility of
effectively sequencing all the numbers in $\Q0$ is primarily
contingent on some knowledge about $\nu_m$, namely the number of the
non-blank entries in the $m^{th}$ row of Table\myref{apptable}

Without pretending to treat thoroughly this topic, we will just
confine ourselves to a few remarks about some simple properties of
the numbers $\nu_m$ (the number of rational numbers present in a row
of Table\myref{apptable} sharing a common irreducible denominator
$m$) and $\sigma_m$ (the sum of the said rational numbers). First of
all it should be said that the normalization condition\refeq{norm}
can be used to find a procedure to progressively calculate the
values of $\nu_m$. For instance, as stated in the
Section\myref{geometric}, when denominators are geometrically
distributed and numerators are conditionally equiprobable, the
distribution of $Q$ is\refeq{geomdistr} and the
normalization\refeq{norm} must be enforced by taking into account
the number $\nu_m$ of the equiprobable numbers sharing the same
irreducible denominator. It is easy to see then that, by taking
$z=1-w$ in\refeq{geomdistr}, the normalization condition\refeq{norm}
becomes
\begin{equation*}
    \sum_{q\in\Q0}\pr{Q=q}=\frac{1-z}{z^2}\sum_{m=1}^\infty\nu_m\sum_{\ell=1}^\infty\frac{z^{m\ell+1}}{m\ell+1}=1
\end{equation*}
namely with a power expansion
\begin{equation}\label{recurr}
   \sum_{m=1}^\infty\nu_m\sum_{\ell=1}^\infty\frac{z^{m\ell+1}}{m\ell+1}=\frac{z^2}{1-z}=\sum_{j=0}^\infty z^{j+2}
\end{equation}
This relation can be used to find the values of $\nu_m$ by equating
the coefficients of the identical powers of $z$: by explicitly
writing indeed the first terms of\refeq{recurr} we find
\begin{align*}
  & \nu_1\left(\frac{z^2}{2}+\frac{z^3}{3}+\frac{z^4}{4}+\frac{z^5}{5}+\ldots\right)
          +\nu_2\left(\frac{z^3}{3}+\frac{z^5}{5}+\frac{z^7}{7}+\frac{z^9}{9}+\ldots\right)\\
  & \qquad\qquad+\nu_3\left(\frac{z^4}{4}+\frac{z^7}{7}+\frac{z^{10}}{10}+\frac{z^{13}}{13}+\ldots\right)
          +\nu_4\left(\frac{z^5}{5}+\frac{z^9}{9}+\frac{z^{13}}{13}+\frac{z^{17}}{17}+\ldots\right)+\ldots\\
  & \qquad=z^2+z^3+z^4+z^5+\ldots
\end{align*}
and hence we progressively have
\begin{align*}
    &^{\nu_1}/_2=1  & \nu_1=2\\
    &^{\nu_1}/_3+\,\!^{\nu_2}/_3=1 & \nu_2=1\\
    &^{\nu_1}/_4+\,\!^{\nu_3}/_4=1 & \nu_3=2\\
    &^{\nu_1}/_5+\,\!^{\nu_2}/_5+\,\!^{\nu_4}/_5=1 & \nu_4=2\\
    &\ldots &\ldots
\end{align*}
and so on, in apparent agreement with the corresponding entries of
the Table\myref{apptable}. It must be added that this procedure can
not be contingent on the specific distribution of $Q$ because
$\nu_m$ is always the same sequence and the normalization
condition\refeq{norm} must hold for every legitimate distribution

We will finally list a few elementary properties of $\nu_m$ and
$\sigma_m$ that can be helpful for every future advance: here
$m=1,2,\ldots$ are the denominators, $n=0,1,\ldots,m$ the numerators
and we call them \emph{accepted} when $^n/_m$ appears in the
Table\myref{apptable}, namely if it is an irreducible fraction:
\begin{enumerate}
    \item \textbf{$\bm{\nu_m\le m-1}$ for $\bm{m\ge2}$}: in our table $n=0$ and $n=m$
are accepted only for $m=1$ so that in every row with $m\ge2$ the
first and last number are always missing; then apparently
$\nu_m=(m+1)-2=m-1$; in particular \textbf{$\bm{\nu_m=m-1}$ only for
$\bm m$ prime number}
    \item \textbf{for $\bm{m\ge3}$, if $\bm{n=k\ge1}$ is accepted, then also $\bm{n=m-k\le m-1}$ is
accepted} because, if $^k/_m$ is irreducible, then also
$^{(m-k)}/_m=1-\,\!^k/_m$ is irreducible, namely the accepted values
always show up in \emph{pairs}; in particular, since $n=1$ is always
accepted, then also $n=m-1$ is always accepted and hence
\textbf{$\bm{\nu_m\ge2}$ for $\bm{m\ge3}$} (the two numbers coincide
for $m=2$, so that $\nu_2=1$)\label{point2}
    \item \textbf{$\bm{\nu_m}$ always is an even number for
$\bm{m\ge3}$} because according to the point\myref{point2} the
accepted numerators $n$ always show up in pairs; moreover \textbf{if
$\bm{m\ge3}$ is even, then $\bm{n=\,\!^m/_2}$ is not accepted}
because for $m=2\ell$ (and $\ell\ge2$) the numerator would be
$n=\,\!^m/_2=\ell$, and $^n/_m=\,\!^\ell/_{2\ell}$ would be a
reducible fraction
    \item \textbf{for $\bm{m\ge3}$ the sum of an accepted pair always is $\bm1$}
because we are adding $^k/_m$ and $^{(m-k)}/_m=1-\,\!^k/_m$; as a
consequence \textbf{the sum of the irreducible fractions sharing a
common denominator $\bm m$ is $\bm{\sigma_m=\,\!^{\nu_m}/_2}$}
because there are $\,^{\nu_m}/_2$ accepted pairs; looking moreover
at the Table\myref{apptable} we see that this last result holds also
for $m=1$ ($\nu_1=2,\,\sigma_1=1$) and $m=2$
($\nu_2=1,\,\sigma_2=\,\!^1/_2$)
\end{enumerate}

%\end{appendix}

%\vspace{10pt}

%\noindent \textbf{Acknowledgements}: The author would like to thank
%a few friends for stimulating comments and suggestions

%\begin{appendix}

%\section{Change of measure}\label{change}

%From the usual textbooks~\cite{pap}, with a slightly adapted
%notation

%\end{appendix}


\begin{thebibliography}{99}

\bibitem{grad}\emph{I.S.\ Gradshteyn and I.M.\ Ryzhik}, \textsc{Table of
Integrals, Series and Products} (Academic Press, Burlington {2007})

%\bibitem{abram} \emph{M.\ Abramowitz and I.\ A.\ Stegun} \textsc{Handbook of Mathematical Functions}
%(Dover Publications,  1968).

%\bibitem{bert} \emph{J.\ Bertrand}, \textsc{Calcul des
%Probabilit\és} (Gauthier-Villars, Paris 1889)

%\bibitem{gned} \emph{B.V. Gnedenko}, \textsc{The Theory of Probability} (MIR, Moscow 1978)

%\bibitem{mathai} \emph{A.M.\ Mathai}, \textsc{An Introduction to Geometrical
%Probability} (Gordon and Breach, Amsterdam 1999)

%\bibitem{pap} \emph{A.\ Papoulis and S.U.\ Pillai}, \textsc{Probability, Random Variables and Stochastic Processes} (McGraw-Hill, Boston 2002)

%\bibitem{czu} \emph{E.\ Czuber}, \textsc{Wahrscheinlichkeitsrechnung und ihre Anwendungen auf Fehlerausgleichung, Statistik und
%Lebensversicherung} (Teubner, Leipzig 1903)

%\bibitem{defin} \emph{B.\ De Finetti}, \textsc{Theory of Probability} (J. Wiley\&Sons, New York 1990)

%\bibitem{shir} \emph{A.N.\ Shiryayev}, \textsc{Probability} (Springer, New York 1984)



%\bibitem{loeve} \emph{M.\ Lo\`eve}, \textsc{Probability Theory I--II} (Springer,
%Berlin, 1977-8).

%\bibitem{jacquet}
%\emph{P.\ Jacquet, W.\ Szpankowski}, IEEE Trans. Inform. Th. 45
%(1999) 1072

%\bibitem{cichon} \emph{J.\ Cicho\'n, Z.\ Go{\l}\c{e}biewski}, DMTCS Proc.
%AQ (2012) 179, 23rd Intern. Meeting on \emph{Probabilistic,
%Combinatorial, and Asymptotic Methods for the Analysis of
%Algorithms} (AofA'12)

%\bibitem{cover} \emph{T.M.\ Cover, J.M.\ Thomas}, \textsc{Elements of Information
%Theory}, (Wiley, Hoboken, New Jersey, 2006)

%\bibitem{betten} \emph{L.M.A.\ Bettencourt, V. Gintautas, M.I.\ Ham},
%Phys.\ Rev.\ Lett. \textbf{100}, 238701 (2008)

%\bibitem{tsallis} \emph{C.\ Tsallis}, J.\ Stat.\ Phys.\ \textbf{52} (1988) 479

%\bibitem{doob}
%\emph{J.L.\ Doob}, \textsc{Stochastic processes} (Wiley, New York
%1953)

%\bibitem{tichonov}
%\emph{A.N.\ Tichonov e A.A.\ Samarskij}, \textsc{Equazioni della
%Fisica Matematica} (Edizioni Mir, Mosca, 1981)

%\bibitem{gardiner} \emph{C.W.\ Gardiner}, \textsc{Handbook of stochastic
%methods} (Springer, Berlin, 1996)


%\bibitem{cufaro09} \emph{N.\ Cufaro Petroni and M.\ Pusterla}, Physica A 388 (2009)
%824.

%\bibitem{nelson}
%\emph{E.\ Nelson}, \textsc{Dynamical theories of Brownian motion}
%(Princeton UP, Princeton 1967) in the new \emph{pdf} version
%available on
%\textsl{web.math.princeton.edu/{\footnotesize$\sim$}nelson/books.html}

%\bibitem{andricuf}
%\emph{A.\ Andrisani and N.\ Cufaro Petroni}, J.\ Math.\ Phys.\ 52,
%113509 (2011)


%\bibitem{feynman}
%\emph{R.\ P.\ Feynman} and \emph{A.\ R.\ Hibbs}, \textsc{Quantum mechanics and path
%integrals} (McGraw--Hill, New York, 1965).\\
%\emph{L.\ S.\ Schulman}, \textsc{Techniques and applications of path integration} (Wiley, New York 1981).\\
%\emph{M.\ Nagasawa}, \textsc{Scr\"odinger equations and diffusion theroy} (Birkh\"auser, Basel
%1993).)

%\bibitem{bohmvigier}
%\emph{D.\ Bohm and J.-P.\ Vigier}, Phys.\ Rev.\ 96 (1954) 208.\\
%\emph{N.\ Cufaro Petroni and F.\ Guerra} Found.\ Phys.\ 25 (1995)
%297.

%\bibitem{guerra}
%\emph{F.\ Guerra}, Phys.\ Rev.\ 77 (1981) 263.\\
%\emph{F.\ Guerra and L.\ Morato}, Phys.\ Rev.\ D 27 (1983)
%1774.\\
%\emph{N.\ Cufaro Petroni and L.\ Morato}, J.\ Phys.\ A 33 (2000) 5833,

%\bibitem{garbaczewski1} \emph{P.\ Garbaczewski, J.\ R.\ Klauder and R.\ Olkiewicz}, Phys.\ Rev.\
%E 51 (1995) 4114.\\
%\emph{P.\ Garbaczewski and R.\ Olkiewicz}, Phys.\ Rev.\ A 51
%(1995) 3445.\\
%\emph{P.\ Garbaczewski and R.\ Olkiewicz}, J.\ Math.\ Phys.\ 40 (1999) 1057.\\
%\emph{P.\ Garbaczewski and R.\ Olkiewicz}, J.\ Math.\ Phys.\ 41 (2000) 6843.\\
%\emph{P.\ Garbaczewski}, Physica A 389 (2010) 936.

%\bibitem{laskin} \emph{N.\ Laskin}, Phys.\ Rev.\ E 62 (2000) 3135.\\
%\emph{N.\ Laskin}, Phys.\ Rev.\ E 66 (2002) 056108.

%\bibitem{applications}
%\emph{S.\ Albeverio, Ph.\ Blanchard and R.\ H\o gh-Krohn},
%Expo.\ Math.\ 4 (1983) 365.\\
%\emph{N.\ Cufaro Petroni, S.\ De Martino, S.\ De Siena and F.\
%Illuminati}, Phys.\ Rev.\  E 63 (2000) 016501.\\
%\emph{N.\ Cufaro Petroni, S.\ De Martino, S.\ De Siena and F.\ Illuminati},
%Phys.\ Rev.\ ST Accel.\ Beams.\ 6 (2003) 034206.\\
%\emph{N.\ Cufaro Petroni, S.\ De Martino, S.\ De Siena and F.\ Illuminati}, Phys.\
%Rev.\ E 72 (2005) 066502.\\
%\emph{N.\ Cufaro Petroni, S.\ De Martino, S.\ De Siena and F.\ Illuminati},
%Nucl.\ Instr.\ Meth.\ A 561 (2006) 237.\\
%\emph{W.\ Paul and J.\ Baschnagel}, \textsc{Stochastic processes: from physics to finance}
%(Springer, Berlin 1999).

%\bibitem{mainardi} \emph{R.\ Gorenflo and F.\ Mainardi} Frac.\ Calc.\
%Appl.\ An. 1 (1998) 167 (reprinted at
%http://www.fracalmo.org/).\\
%\emph{R.\ Gorenflo and F.\ Mainardi} Arch.\ Mech.\ 50 (1998) 377 (reprinted at
%http://www.fracalmo.org/).

%\bibitem{cufaro10arXiv} \emph{N.\ Cufaro Petroni}, (2010) arXiv:1006.1833v1 [math.PR].

%\bibitem{sato} \emph{K.\ Sato}:  \textsc{L\'evy processes and infinitely
%divisible distributions} (Cambridge U.P., Cambridge, 1999)

%\bibitem{applebaum} \emph{D.\ Applebaum}: \textsc{L\évy processes and Stochastic Calculus}
%(Cambridge U.P., Cambridge, 2004).

%\bibitem{cufaro08} \emph{N.\ Cufaro Petroni}, Physica A 387 (2008)
%1875.

%\bibitem{cufaro07} \emph{N.\ Cufaro Petroni}, J.\ Phys.\ A 40 (2007)
%2227.

%\bibitem{BergVignat} \emph{C.\ Berg and C.\ Vignat}, J.\ Phys.\ A 41 (2008) 265004

%\bibitem{HeideLeonenko} \emph{C.\ C.\ Heyde and N.\ N.\ Leonenko}, Adv.\ Appl.\ Prob.\
%37 (2005) 342.

%\bibitem{Feller2} \emph{W.\ Feller}, \textsc{An introduction to probability theory and its
%applications II} (Wiley\& Sons, New York, 1971).

%\bibitem{abramowitz} \emph{M.\ Abramowitz and I.\ A.\ Stegun} \textsc{Handbook of Mathematical Functions}
%(Dover Publications,  1968).

%\bibitem{balakrishnan} \emph{N.\ Balakrishnan and V.\ B.\ Nevzorov},
%\textsc{A primer on statistical distributions} (Wiley\& Sons,
%Hoboken, 2003).

%\bibitem{rudin} \emph{W,\ Rudin}, \textsc{Functional analysis} (McGraw-Hill, New York, 1973)

%\bibitem{stroock} \emph{D.\ W.\ Stroock}, \textsc{An introduction to partial
%differential equations for probabilists} (Cambridge U.P., Cambridge,
%2008).

%\bibitem{varadhan} \emph{D.\ W.\ Stroock and S.\ R.\ S.\ Varadhan},
%\textsc{Multidimensional diffusion processes}, (Springer, Berlin,
%2006).

%\bibitem{madelung} \emph{E.\ Madelung}, Z.\ Physik 40 (1926) 332.

%\bibitem{vivoli} \emph{A.\ Vivoli, C.\ Benedetti and G.\ Turchetti},
%Nucl.\ Instr.\ Meth.\ A 561 (2006) 320.

%\bibitem{paul} \emph{W.\ Paul and J.\ Baschnagel} \textsc{Stochastic processes: from physics
%to finance} (Springer, Berlin 1999).

%\bibitem{protter} \emph{Ph.\ E.\ Protter}, \textsc{Stochastic integration and
%differential equations} (Springer, Berlin 2005).\\
%\emph{O.\ E.\ Barndorff--Nielsen et al} ed's, \textsc{L\'evy
%processes, Theory and applications} (Birkh\"auser, Boston 2001).\\
%\emph{N.\ Jacob}, \textsc{Pseudo--differential operators and Markov
%processes} vol I--III (Imperial College Press, London 2001--05).


\end{thebibliography}
\end{document}